\documentclass[11pt]{article}
\usepackage{amsmath,amsthm,amssymb,amscd,fancybox,boxedminipage,shadow,ifthen}
\usepackage{epsfig,subfigure}
\usepackage{times}

\newcommand{\C}{{\mathbb C}}
\newcommand{\half}{{\textstyle{\frac{1}{2}}}}
\newcommand{\Li}{{\rm Li}}
\newcommand{\li}{{\rm li}}
\newcommand{\R}{{\mathbb R}}
\newcommand{\sgn}{{\rm sgn}}
\newcommand{\Z}{{\mathbb Z}}
\newcommand{\g}{{\gamma}}

\newcommand\cH{\mathcal H}

\newcommand\cF{\mathcal F}

\newcommand\cS{\mathcal{S}}

\newcommand\hvarphi{\hat\varphi}

\newcommand\hpsi{\hat\psi}

\newcommand\hf{\hat{f}}

\newcommand\datver[1]{\def\datverp%
 {\par\boxed{\boxed{\text{#1; Run: \today}}}}}
\newcommand\ha{\frac12}

\renewcommand\Re{\operatorname{Re}}

\newcommand\PV{\operatorname{PV}}

\newcommand\bbC{\mathbb C}

\newcommand\bbN{\mathbb N}

\newcommand\bbR{\mathbb R}

\newcommand\pa{\partial}

\newcommand\ch[2]{\left(\begin{matrix} #1\\#2\end{matrix}\right)}

\newcommand\CI{{\mathcal C}^{\infty}}

\newtheorem{theorem}{Theorem}
\newtheorem{proposition}{Proposition}
\newtheorem{corollary}{Corollary}

\theoremstyle{definition}
\newtheorem{definition}{Definition}

\theoremstyle{remark}
\newtheorem{remark}{Remark}

\datver{Version 0.1, September 18, 2006}
\begin{document}

\title{Tempering the Polylogarithm}

\author{Charles L. Epstein and Jack Morava
\footnote{Research of both authors partially supported by DARPA under the
FUNBIO program.}\\ University of Pennsylvania and\\ Johns Hopkins University}

\maketitle 

\begin{abstract} We show that the function $\Li_s(e^x)$ extends to
an entire function of the complex variable s, taking values in tempered
distributions (in $x$) on the whole real line.  As a corollary the classical
polylogarithm extends to an entire function (of $s$) taking values in
distributions on the positive real axis.  We then identify the singularities of
$\Li_s(e^x)$ in terms of distributional powers of x, which leads to a simple proof
of the smoothness of the `modified' polylogarithm of Bloch, Ramakrishnan,
Wojtkowiak, Zagier, and others.
\end{abstract}

\section{The polylogarithm}
The power series
\[
\Li_s(z) = \sum_{n \geq 1} \frac{z^n}{n^s}
\]
converges when $|z| <1$, defining the classical polylogarithm function, equal to 
$ - \log (1 - z)$ when $s=1$. In general, the behavior of these functions at 
$z = 1$ is complicated: it is known, for example, that $\Li_s$ has an analytic 
continuation to the cut plane $\C - [1,\infty)$. The fact that the `modified' 
polylogarithm of Bloch, Ramakrishnan, Wigner, Wojtkowiak, Zagier, \dots, defined 
as the real, or imaginary part of 
\[
\sum_{n-1 \geq k \geq 0} \frac{(-2)^kB_k}{k!} \; \log^k|z| \cdot \Li_{n-k}(z),
\]
according to whether $n \in \Z$ is even or odd, see~\cite{EGK}, extends to
define a smooth function on the whole complex plane gives some idea of
the complexity of its branch-point behavior.  

The closely related series 
\[
\li_s(x) = \sum_{n \geq 1} \frac{e^{nx}}{n^s} = \Li_s(e^x) 
\]
converges when $\Re s \geq 0$ and $\Re x < 0.$ This note is concerned with the
function, or rather the tempered distribution, it defines upon restriction to
the real axis. In~\cite{Morava}, the second author considered the
polylogarithm on the unit circle; our recent interest in its behavior on the
real line comes from problems in statistical mechanics, see~\cite{morava2}. 

We would like to thank Don Zagier for pointing out that, when $\Re s < 1$, the
methods of~\cite{Zagier} can be used to establish that
\[
\li_s(x) =  \Gamma(1-s)(-x)^{s-1} + \sum_{k \geq 0} \zeta(s-k) \frac{x^k}{k!} 
\]
(if $|x| < 2 \pi$ as well, which ensures convergence of the series on the right).
Our starting point is the simpler observation that, when $\Re s \geq 1,$ and $x<0,$
\[
\pa_x\li_s(x) = \li_{s-1}(x).
\]

It is useful to establish some notation. Recall that convolution of a smooth 
function supported on the positive real axis with the distribution-valued
divided power
\[
\gamma^{s}_+(x) = \frac{x_+^{s-1}}{\Gamma(s)}
\]
defines an entire function of the complex variable $s,$ see~\cite{GelShi}[Ch. I
\S 3.5], representing fractional differentiation of order $-s,$ see \cite{GelShi}[Ch.
I \S 5.5].

\begin{proposition}\label{prop1.0} If $\Re s > 0$ and $0 >  x \in \R$, then 
\[
\li_s(x) = (\gamma^{s}_+ \; * \; \li_0)(x) 
\]
as functions holomorphic in the right half $s$-plane, with values in smooth functions 
of $x <0$; where 
\[
\li_0(x) = \frac{1}{e^{-x} - 1} = -\frac{1}{x} \sum_{n \geq 0} B_n \frac{x^n}{n!} \;.
\]
\end{proposition}

\begin{proof} The classical integral representation
\[
\sum_{n \geq 0} \frac{z^n}{(n+a)^s} = \frac{1}{\Gamma(s)} \int_0^\infty \frac{t^{s-1}
e^{-at}}{1 - ze^{-t}} \; dt
\]
for the `Lerch transcendent', valid for $\Re a > 0, \; \Re s > 0$, and $|z| < 1$, is 
obtained by expanding the denominator in the integral as a power series. Take $a=1, 
\; z = e^x$ with $x < 0$, and multiply both sides by $z$ to get
\[
\li_s(x) = \int_{-\infty}^\infty \frac{t_+^{s-1}}{\Gamma(s)} \; \frac{e^{-(t-x)}}
{1 - e^{-(t-x)}} \; dt = (\gamma^{s}_+ \; * \; \li_0)(x).
\]
\end{proof}

In our case $\gamma^{s}_+$ has at worst polynomial growth on the right
half-line,  $\li_0(x) \to 0$ exponentially as $x \to -\infty,$ and 
$\li_0(x) +1 \to 0$ exponentially as $x \to \infty.$  Near $x=0$, however, 
\[
\li_0 \equiv - x^{-1} \; {\rm mod \; smooth \; functions}
\]
is singular, and the behavior of the convolution there is quite interesting. 


\section{Extension of $\li_s$ as a tempered distribution}
Much of the material in this section is standard and can be found
in~\cite{Hormander2}. We let $\cS(\bbR)$ denote the Schwartz class functions on $\bbR$
with topology defined by the semi-norms
\begin{equation}
\|f\|_N=\max_{0\leq j\leq N}\sup_{x\in\bbR} (1+|x|)^N|\pa_x^jf(x)|, N\in\bbN.
\end{equation}
The dual space $\cS'(\bbR)$ is the space of tempered distributions, with the
usual weak topology. This implies that $l\in\cS'(\bbR)$ if and only if there is
an $N$ so that
\begin{equation}
|l(f)|\leq C_N\|f\|_N.
\end{equation}
Hence a given tempered distribution is always of finite order, and growth at
infinity, and therefore can be continuously extended to a much larger
space. Schwartz space is identified with a subset of $\cS'(\bbR)$ by $f\mapsto l_f:$
\begin{equation}
l_{f}(\varphi)=\int \varphi(x) f(x)dx.
\end{equation}
It is a very useful and important fact that $\cS(\bbR)$ is a \emph{dense}
subspace.

Recall that  if $\varphi$ is a tempered distribution,  then the Fourier
transform $\hvarphi$ is the tempered distribution \emph{defined} by duality, with:
\begin{equation}
\langle \hvarphi,f\rangle=\langle\varphi,\hf\rangle \text{ for all }f\in\cS(\bbR).
\end{equation}
This defines a tempered distribution, because the Fourier transform is an
isomorphism of $\cS(\bbR)$ to itself.  In general the product of two
distributions is not defined, and because $\widehat{\varphi
*\psi}=\hvarphi\hpsi,$ this implies that the convolution of two tempered
distribution is not always defined. Of course it may be defined for a given
pair.

When it is defined, similar considerations are used to define the convolution of two
distributions. We start with the case
$\varphi,\psi\in\cS(\bbR)$ and observe that, for $f\in\cS(\bbR)$ we have:
\begin{equation}
\begin{split}
\langle \varphi *\psi,f\rangle&=\int\left[\int\varphi(x)\psi(t-x)dx\right] f(t)dt\\
&=\int\varphi(x)\left[\int\psi(y) f(x+y)dy\right] dx
\end{split}
\end{equation}
We changed the order of integrations and variables, with $t-x=y$ to go from the
first line to the second. Because $\cS(\bbR)\subset \cS'(\bbR)$ is dense, this
relation defines $\psi*\varphi,$ whenever this convolution defines a tempered
distribution. For $x\in\bbR$ we define the operation
$\tau_x:\cS(\bbR)\to\cS(\bbR)$ by
\begin{equation}
\tau_x f(y)=f(x+y).
\end{equation}
The convolution of $\varphi$ and $\psi$ defines a tempered distribution
precisely when 
\begin{enumerate}
\item $x\mapsto \langle\varphi,\tau_x f\rangle$ is in the domain of $\psi.$
\item There is an $N$ and a $C_N>$ so that
\begin{equation}
|\langle \psi,\langle\varphi,\tau_x f\rangle\rangle|\leq C_N\|f\|_{N}.
\end{equation}
\end{enumerate}
In this case the convolution $\varphi *\psi$ is defined by
\begin{equation}
\langle\varphi *\psi,f\rangle=\langle \varphi(x),\langle\psi,\tau_x
f\rangle\rangle.
\end{equation}

We can now  prove that the tempered distributions 
\begin{equation}
\gamma_+^s(x)=\frac{x_+^{s-1}}{\Gamma(s)},\text{ and }\li_0=\PV\frac{1}{e^{-x}-1},
\end{equation} 
can be convolved to give an entire family of tempered distributions. Observe
that, the discussion above implies that if $\gamma_+^s*\li_0$ makes sense as a
distribution, then for all Schwartz class functions $f,$ we must have
\begin{equation}
\langle
\gamma_+^s*\li_0,f\rangle=\PV\int\limits_{-\infty}^{\infty}\frac{1}{e^{-t}-1}
\langle \gamma_+^s(\cdot), f(\cdot+t)\rangle dt.
\end{equation}

The crux of the matter is therefore to analyze
\begin{equation}
F_s(t)=\langle \gamma_+^s(\cdot),f(\cdot+t)\rangle,
\end{equation}
as a function of $(t,s).$ If $\Re s>1,$ then
$F_s(t)$ is given by an absolutely convergent integral and for any $k\in\bbN,$
we can integrate by parts to obtain:
\begin{equation}
F_s(t)=\frac{(-1)^k}{\Gamma(s+k)}\int\limits_{0}^{\infty}x^{s+k-1}f^{[k]}(x+t)dx.
\label{12}
\end{equation}
The right hand is an analytic function in $-k<\Re(s)<k$ with values in the
space of functions, $\cF^+_k,$ which we now define:
\begin{definition} A function $f\in\CI(\bbR)$ belongs to $\cF^+_k$ if 
\begin{enumerate}
\item $f(x)$ and all its derivatives are rapidly decreasing as $x$ tends to
  $+\infty.$
\item For all $j$ 
\begin{equation}
\limsup_{x\to-\infty}\left|\frac{\pa_x^jf(x)}{(1+|x|)^{k}}\right|<\infty.
\end{equation}
\end{enumerate}
\end{definition}
Briefly, $f\in\cF^+_k$ if $f$ is smooth, in Schwartz class at $+\infty$ and,
of tempered growth at $-\infty.$  The topology on $\cF^+_k$ is defined by the
semi-norms:
\begin{equation}
|F|_{k,l}=\sup_{x>0}[(1+|x|)^l\max_{0\leq j\leq l}|\pa_x^j F(x)|]+
\sup_{x\leq 0}[(1+|x|)^{-k}\max_{0\leq j\leq l}|\pa_x^j F(x)|].
\end{equation}

The statement  that $F_s$ is an analytic function from $-k<\Re s<k$ to $\cF^+_k$ is now
a simple consequence of formula~\eqref{12}. Moreover, $f\mapsto F_s$ is
continuous as a mapping from $\cS(\bbR)$ to $\cF^+_k,$ that is for each $l$
there is an $N_l$ and a $C_{s,l}$ so that
\begin{equation}
|F_s|_{k,l}\leq C_{s,l}\|f\|_{N_l}.
\end{equation}
The constants $C_{s,l}$ are locally uniformly bounded in $|\Re s|<k.$ 

To
complete our discussion of $\gamma_+^s *\li_0$ it remains only to show that
$\cF^+_k$ is in the domain of $\li_0$ for all $k\in\bbN.$ To that end we
choose an even function $\chi\in\CI_c((-1,1)),$ which equals $1$ in the
interval $[-\ha,\ha].$ For any such function we have
\begin{equation}
\langle\li_0,f\rangle=\int\limits_{-\infty}^{\infty}\frac{(1-\chi(t))f(t)dt}{e^{-t}-1}+
\lim_{\epsilon\to
  0^+}\int\limits_{\epsilon<|t|}\frac{\chi(t)f(t)dt}{e^{-t}-1}.
\label{14}
\end{equation}
Observe that
\begin{equation}
\frac{1}{e^{-t}-1}=\begin{cases} -1+O(e^{-t})&\text{ as }t\to+\infty\\
O(e^{-|t|})&\text{ as }t\to -\infty
\end{cases}
\end{equation}
>From these estimates it follows that the first term on the r.h.s. of~\eqref{14}
is clearly a continuous functional on $\cF^+_k,$ for every $k.$

We consider the principal value part. An elementary calculation shows that, for
a differentiable $f$ we have:
\begin{equation}
\lim_{\epsilon\to
  0^+}\int\limits_{\epsilon<|t|}\frac{\chi(t)f(t)dt}{e^{-t}-1}=
\int\chi(t)f(t)\left[\frac{1}{e^{-t}-1}+\frac{1}{t}\right]dt-
\int\chi(t)\frac{f(t)-f(0)}{t}dt.
\label{18}
\end{equation}
The sum $(e^{-t}-1)^{-1}+t^{-1}$ is a smooth function, hence the right hand
side again clearly defines a continuous functional on $\cF^+_k,$ for every $k.$
This proves the theorem
\begin{theorem} The family $s\mapsto \gamma_+^s *\li_0$ is an entire family of tempered
  distributions.
\end{theorem}
In the sequel we use $\li_s$ to denote the distribution $ \gamma_+^s *\li_0.$
While it is not immediately clear, we will show that the two, a priori
different, definitions of  $\li_0$ do coincide.

The theorem and the well known functional equation,
$\pa_x\gamma_+^s=\gamma_+^{s-1},$ show that the functional equation satisfied
by $\li_s(x)$ in $x<0$ extends to the whole real line.
\begin{corollary}\label{cor1.0} In the sense of distributions, $\pa_x\li_s=\li_{s-1},$ for all
  $s\in\bbC.$
\end{corollary}

\section{The singularities of $\li_s(x).$}

We now turn to a consideration of the singularities of the distribution
$\li_s(x),$ as a function of $x.$ In the previous section we obtained the formula:
\begin{equation}
\langle\gamma_+^s*\li_0,f\rangle= 
\int\limits_{-\infty}^{\infty}\frac{(1-\chi(t))F_s(t)dt}{e^{-t}-1}+
\lim_{\epsilon\to
  0^+}\int\limits_{\epsilon<|t|}\frac{\chi(t)F_s(t)dt}{e^{-t}-1}.
\label{24}
\end{equation}
Using formula~\eqref{12} it is straightforward to
show that 
\begin{equation}
\int\limits_{-\infty}^{\infty}\frac{(1-\chi(t))F_s(t)dt}{e^{-t}-1}=
\frac{1}{\Gamma(s+k)}
\int\limits_{-\infty}^{\infty}\int\limits_{0}^{\infty}
\frac{(1-\chi(x-y))y^{s+k-1}dy}{e^{y-x}-1}f^{[k]}(x)dx.
\end{equation}
Hence there is an analytic family of smooth, tempered \emph{functions} $G_s(x)$ so
that the first term on the right hand side of~\eqref{24} has a representation
as
\begin{equation}
\int\limits_{-\infty}^{\infty}\frac{(1-\chi(t))F_s(t)dt}{e^{-t}-1}=
\langle G_s,f\rangle.
\end{equation}
Thus the singularities of $\li_s(x)$ are  in the second term in~\eqref{24}. 

Using~\eqref{18} we can show that the first term contributes a smooth
term and therefore the distribution, $g_s$ defined by
\begin{equation}
\begin{split}
\langle g_s,f\rangle&=-\int\chi(t)\frac{F_s(t)-F_s(0)}{t}dt\\
&=-\PV\int\chi(t)\frac{F_s(t)dt}{t}.
\end{split}
\label{23}
\end{equation}
has the same singularity as $\li_s.$  

Given a choice of smooth, even cutoff function $\chi$ we define the operator on
Schwartz class functions:
\begin{equation}
\cH_\chi f=-\PV\int\chi(t)\frac{f(x-t)dt}{t}.
\end{equation}
The distribution $h_0=\PV\left[\frac{\chi(t)}{t}\right]$ is compactly supported and
therefore has a smooth Fourier transform, on the other hand
$h_1=\left[\frac{1-\chi(t)}{t}\right]$ is smooth and belongs to $L^2,$ hence its
Fourier transform is rapidly decreasing. The Fourier transform of
$\PV\left[\frac{1}{t}\right]$ is well known to be $-\pi i\sgn\xi.$ This shows that
the Fourier transform of $h_0$ is smooth and rapidly approaches $-\pi
i(\pm 1),$ as $\xi\to\pm\infty.$ Thus $\cH_\chi$ maps $\cS$ to
itself and therefore $\cH_\chi$ extends as a map from $\cS'$ to itself. An
elementary calculation shows that, for $\Re s>-k$ we have
\begin{equation}
\langle g_s,f\rangle=\langle\cH_\chi \gamma_+^{s+k},(-1)^kf^{[k]}\rangle.
\label{255}
\end{equation}
The singularity of $\li_s(x)$ at $x=0$ agrees with that of 
\begin{equation}
g_s=\pa_x^k \cH_\chi \gamma_+^{s+k}=\cH_{\chi}\gamma_+^s.
\end{equation}

It would be tempting to say that this agrees with the singularity of the
Hilbert transform of $\gamma_+^{s},$ but for the fact that the Hilbert
transform does not preserve Schwartz space, and hence does not have an extension to
tempered distributions. Some care is required to compute the singular part of
of $g_s.$ We make extensive usage of the fact that the Fourier transform of a
compactly supported distribution is smooth.

Notice that in~\eqref{255} the power $s+k$ is positive, this,
coupled with the functional equation
\begin{equation}
\pa_x \gamma_+^s=\gamma_+^{s-1}.
\label{256}
\end{equation}
facilitate the computations which follow. We first compute the Fourier
transform of $\gamma_+^s$
\begin{proposition}
The Fourier transform of the tempered distribution $\gamma_+^s$ is
$e^{-i\frac{\pi s}{2}}\eta_{-}^{-s},$ where
\begin{equation}
\eta_{\pm}^{s}(\xi)=\lim_{\epsilon\downarrow 0}(\xi\pm i\epsilon)^s.
\end{equation}
The complex power $z\mapsto z^s$ is taken to be real on the positive real axis
and analytic in $\bbC\setminus (-\infty,0].$
\end{proposition}
\begin{proof}
If $s>0,$ then we can compute the Fourier transform using
\begin {equation}
\begin{split}
\widehat{\gamma_+^s}(\xi)&=\lim_{\epsilon\downarrow 0}\frac{1}{\Gamma(s)}
\int\limits_{0}^\infty x^{s-1}e^{-\epsilon x}e^{-ix\xi}dx\\
&=\lim_{\epsilon\downarrow 0}\frac{1}{\Gamma(s)}
\int\limits_{0}^\infty x^{s-1}e^{- x(\epsilon+i\xi)}dx.
\end{split}
\end{equation}
Given our choice of branch for $z^{s-1},$ the last integral can be regarded as
a contour integral along the ray $\{x(\epsilon+i\xi):\: x>0\},$ which lies in the right
half plane. The conclusion follows from an elementary contour deformation
argument. For $s\leq 0$ we use the functional equation to conclude that
\begin{equation}
\widehat{\gamma_+^{s-k}}=(i\xi)^k\widehat{\gamma_+^{s}},
\end{equation}
which is easily seen to extend our formula for $\widehat{\gamma_+^{s}}$ to the
whole complex plane.
\end{proof}

The Fourier transform of $\cH_\chi\gamma_+^s$ is $-h_0(\xi)e^{-i\frac{\pi
    s}{2}}\eta_-^{-s}.$ Let $\psi(\xi)$ be a smooth, even, non-negative function which
    vanishes in a neighborhood of $0,$ and equals $1$ for $|\xi|>1,$ then
\begin{equation}
\cH_{\chi}\gamma_+^s=-e^{-i\frac{\pi
    s}{2}}\left[\cF^{-1}((1-\psi)h_0(\xi)\eta_-^{-s})+
\cF^{-1}(\psi h_0(\xi)\eta_-^{-s})\right]
\end{equation}
The distribution $(1-\psi)h_0(\xi)\eta_-^{-s}$ is compactly supported, hence
its inverse Fourier transform is a smooth function. The singularity of $g_s$ is
therefore equal to that of
\begin{equation}
g_{s0}=-e^{-i\frac{\pi s}{2}}\cF^{-1}(\psi h_0(\xi)\eta_-^{-s}).
\end{equation}
>From the remarks above, it is clear that the difference $\psi
h_0(\xi)-\psi(-i\pi\sgn\xi)$ is a smooth rapidly decreasing function, and
therefore the singularity of $g_{s0}$ equals that of
\begin{equation}
g_{s1}=i\pi e^{-i\frac{\pi s}{2}}\cF^{-1}(\psi \sgn\xi\eta_-^{-s}).
\end{equation}

If we  let $\psi_+(\xi)=\chi_{[0,\infty)}(\xi)\psi(\xi),$ then we can write
  $g_{s1}$ as
\begin{equation}
g_{s1}=i\pi\cF^{-1}\left[\frac{e^{-i\frac{\pi s}{2}}\psi_+(\xi)}{\xi^s}-
\frac{e^{i\frac{\pi s}{2}}\psi_+(-\xi)}{|\xi|^s}\right]
\end{equation}
For $\Re s >1,$ we see that
\begin{equation}
g_{s1}=\frac{i}{2}\bigg[e^{-i\frac{\pi s}{2}}\int\limits_{0}^{\infty}
\frac{\psi_+(\xi)e^{ix\xi}d\xi}{\xi^s}-
e^{i\frac{\pi s}{2}}\int\limits_{0}^{\infty}
\frac{\psi_+(\xi)e^{-ix\xi}d\xi}{\xi^s}\bigg].
\end{equation}
Finally, we see that, for all $s\in\bbC,$ the distributions
$\Gamma(1-s)\gamma_+^{1-s}-\frac{\psi_+(\xi)}{\xi^s}$ are compactly supported
and therefore $g_{s1}$ has the same singularity as
\begin{equation}
g_{s2}=
i\pi\Gamma(1-s)\left[e^{-i\frac{\pi s}{2}}\cF^{-1}(\gamma_+^{1-s})(x)-
e^{i\frac{\pi s}{2}}\cF^{-1}(\gamma_+^{1-s})(-x)\right].
\label{38}
\end{equation}

Using calculations similar to the proof of Proposition~\label{prop1} to
evaluate the right hand side of~\eqref{38} we obtain:
\begin{proposition} The function $\li_s$, taking values in smooth functions on the 
negative real line when $\Re s \geq 0$, extends to an entire function of $s$,
taking values in tempered distributions on the whole real line, satisfying the
congruence
\begin{equation}
\begin{split}
\li_s(x) &\equiv -\frac{\Gamma(1-s)}{2}\left[e^{-i\pi s}\eta_+^{s-1}+e^{i\pi
    s}\eta_-^{s-1}
\right]\\
&=-\frac{\Gamma(1-s)}{2}\left[e^{-i\pi s}(x+i0)^{s-1}+e^{i\pi
    s}(x-i0)^{s-1}\right]
\end{split}
\label{big1}
\end{equation}
(modulo meromorphic functions with smooth coefficients).
\end{proposition}
\noindent
Note that when $x < 0$ this accounts precisely for Zagier's correction. \bigskip

\noindent
As follows easily from the computations used in the proof of
Corollary~\ref{cor1}, the value of this distribution at $s=1$ is
\[
\li_1(x) = - \log|1 - e^x|.
\]
Note, this is an identity, not a congruence. Using the functional equation from
Corollary~\ref{cor1.0} and this equation, we can show, as asserted above, that
our notation is consistent. For $x\neq 0,$
\begin{equation}
\pa_x\li_1(x)=\frac{1}{e^{-x}-1}=\li_0(x),
\end{equation}
which implies that, the distribution $\li_s\big|_{s=0},$ given by analytic continuation,
agrees with $\PV(e^{-x}-1)^{-1}.$

The families of distributions $\eta^{s}_{\pm}$ are entire, hence the expression
on the right hand side of equation~\eqref{big1} is holomorphic for $\Re s < 0$,
but has simple poles at positive integers; its residue at $s = n$, however, is
(up to sign) the integral divided power $x^{n-1}/(n-1)!$,
which is smooth. This shows that, while $\li_s$ is itself an entire family of
distributions, its separation into regular and singular parts cannot be done
holomorphically.  Near $s=1$, for example, this implies that
\[
\li_s(x) \equiv \frac{1}{s-1} + \dots \;,
\]
which is related to a similar property of $\zeta(s)$. More generally, \bigskip

\noindent
\begin{corollary}\label{cor1} For a non-negative integer $n$,
\[
\li_{-n}(x) \equiv (-1)^{n-1} n! x^{-n-1}
\]
modulo smooth functions, while for positive integers $n$,
\[
\li_n(x) \equiv \frac{-x^{n-1}\log|x|}{(n-1)!} \;.
\]
\end{corollary}
\noindent
\begin{proof} The first assertion is a consequence of ~\cite{GelShi}[Ch. I \S 3.7], 
but we give  elementary proofs of both formul{\ae}.  These statements follow
from the equation~\eqref{23} for the singular part of $g_s.$ For a
positive integer, we see that
\begin{equation}
\begin{split}
\langle g_n,f\rangle&=-\frac{1}{\Gamma(n)}\sum\limits_{j=0}^{n-1}\ch{n-1}{j}
\PV\int\limits_{-\infty}^{\infty}\frac{\chi(t)(-t)^{j}dt}{t}
\int\limits_{t}^{\infty}y^{n-1-j}f(y)dy\\
&\equiv \frac{-1}{\Gamma(n)}\PV\int\limits_{-\infty}^{\infty}\frac{\chi(t)dt}{t}
\int\limits_{t}^{\infty}y^{n-1}f(y)dy
\end{split}
\end{equation}
For $t>0$ define
\begin{equation}
w(t)=-\int\limits_{t}^{\infty}\frac{\chi(s)ds}{s}\equiv \chi(t)\log t.
\end{equation}
Integrating by parts shows that
\begin{equation}
\begin{split}
\langle g_n,f\rangle&\equiv-\frac{1}{\Gamma(n)}\lim_{\epsilon\downarrow 0}
\left[\int\limits_{-\infty}^{\infty}w(|t|)t^{n-1}f(t)dt+
w(\epsilon)\int\limits_{-\epsilon}^{\epsilon}y^{n-1}f(y)dy\right]\\
&\equiv\int\limits_{-\infty}^{\infty}\frac{-t^{n-1}\log|t|}{(n-1)!}f(t)dt
\end{split}
\end{equation}
as asserted. 

For $-n\leq 0$ we use equation~\eqref{12}, with $k=1+n,$ and~\eqref{23} to
conclude that
\begin{equation}
\langle g_{-n},f\rangle=-\PV\int\frac{\chi(t)(-1)^nf^{[n]}(t)dt}{t}.
\end{equation}
Thus
\begin{equation}
g_{-n}\equiv -\pa_t^n\PV\frac{\chi(t)}{t}.
\end{equation}
The right hand side is equivalent to a homogeneous extension of
$(-1)^nn!t^{-1-n},$ as asserted.
\end{proof}
\noindent
Note that if $n > 2$,
\[
(\g_+^{n-1}(x) \log |x|)' \equiv \g_+^{n-2}(x) \log |x|
\]
modulo smooth functions, and that
\[
\li_1(x) = - \log|1 - e^x| \equiv - \log|x|.
\] 
The Taylor coefficients at zero of the smooth difference can be read off from 
the power series formula
\[
\log(\frac{x}{1 - e^x}) = -\sum_{k \geq 1} \frac{1}{k} \Bigl(\sum_{n \geq 1}(-1)^{n-1} B_n 
\frac{x^n}{n!} \Bigr)^k \;.
\] 

The modified polylogarithm, as defined in the introduction, is 
(essentially) the real or imaginary part, depending on whether $n$ is even or odd, 
of a certain combination of products of classical and polylogarithms, with 
coefficients as above. 
\begin{corollary} When $n$ is even, the distribution analog
\[
\lambda i_n(x) = \sum_{n-1 \geq k \geq 0} B_k \frac{(-2x)^k}{k!} \cdot \li_{n-k}(x) 
\]
of the `modified' polylogarithm is a  smooth function of $x$.
\end{corollary}
\noindent
\begin{remark} The corollary asserts that if we regard the polylogarithm
as a distribution, then this combination is in fact smooth when $n$ is even (which
recovers the smoothness of the classical `modified' polylogarithm). When $n$ is 
odd, the argument below shows that the singular part of the corresponding combination 
is real, so the imaginary part is again smooth. Thus our analysis of the
singularities of the distribution extension of the polylogarithm is precise enough
to recover the smoothness of the classical modified polylogarithm in all
cases.
\end{remark}

\noindent
\begin{proof}
We claim that Corollary~\ref{cor1} implies that
\[
\lambda i_n(x) \equiv - (-2)^{n-1}B^{n-1}(-\half) \cdot \gamma_+^{n-1}(x) \log |x| \;,
\]
where the Bernoulli polynomials are defined by the generating function
\[
\sum_{n \geq 0} B^n(q) \frac{t^n}{n!} = \frac{te^{(q + 1)t}}{e^t - 1}
\]
see~\cite{Iwasawa}[\S 2]. They satisfy
\[
B^n(q) = \sum_{n \geq k \geq 0} \ch{n}{k} B_k q^{n-k} \;.
\]
Substituting the expression for the singular part of $\li_{n-k}$ implies that
\[
\lambda i_n(x) \equiv - \sum_{n-1 \geq k \geq 0} \frac {(-2)^k B_k}{k! (n-k-1)!} 
\cdot x^{n-1} \log |x| \;,
\]
and rewriting $(-2)^k$ as $(-\half)^{1-n}(-\half)^{n-k-1}$ yields the assertion. 
It follows from the definition that 
\[
\sum_{n \geq 0} B^n(-\half) \frac{t^n}{n!} = \frac{t}{2 \sinh \half t} 
\]
is an even function, so $B^{\rm odd}(-\half) = 0$.
\end{proof}
\noindent
It is obvious, on the other hand, that $B^{\rm even}(-\half) \in \R$.

As a final corollary we observe that:
\begin{corollary} The function $s \mapsto \Li_s$ extends to an entire function 
of $s$, taking values in distributions on the positive real axis, defined on 
compactly-supported test functions by
\[
f \mapsto \int_0^\infty \Li_s(t)f(t) \; dt = \int_{-\infty}^\infty \li_s(x) \cdot 
e^x f(e^x) \; dx \;.
\]
\end{corollary}
\begin{remark} The domain of $\Li_s$ contains the set of functions such that
  $e^xf(e^x)\in\cS(\bbR).$ For example a smooth function on $[0,\infty)$ such
  that, for every $k\in\bbN,$ there is a $C_k$ satisfying
\begin{equation}
\sup_{y\geq 0}|(1+y^2)\pa_y^kf(y)|\leq C_k,
\end{equation}
belongs to the domain of $\Li_s.$
\end{remark}

\end{document}